\documentclass[11pt,reqno]{amsart}
\usepackage{amsthm}
\usepackage{amssymb}
\usepackage{graphicx}
\usepackage{latexsym}
\usepackage{multicol}
\usepackage{verbatim,enumerate}
\usepackage{accents}
\usepackage[usenames]{color}
\usepackage[colorlinks=true, linkcolor=blue, citecolor=blue, urlcolor=blue]{hyperref}
\usepackage{hyperref}
\usepackage{amsmath, amscd}

\theoremstyle{definition}
\newtheorem*{cor}{Corollary}
\newtheorem*{lem}{Lemma}
\newtheorem*{prop}{Proposition}

\newtheorem{thm}{Theorem}

\theoremstyle{definition}
\newtheorem*{defn}{Definition}

\newtheorem*{example}{Example}

\newtheorem*{rem}{Remark}

\newtheorem*{theom}{Theorem}

\newcommand{\nc}{\newcommand}

\newcounter{cnt}
\numberwithin{equation}{section}

\makeatletter
\def\section{\def\@secnumfont{\mdseries}\@startsection{section}{1}%
  \z@{.7\linespacing\@plus\linespacing}{.5\linespacing}%
  {\normalfont\scshape\centering}}
\def\subsection{\def\@secnumfont{\bfseries}\@startsection{subsection}{2}%
  {\parindent}{.5\linespacing\@plus.7\linespacing}{-.5em}%
  {\normalfont\bfseries}}
\makeatother

 \nc{\Hom}{\operatorname{Hom}}
  \nc{\mode}{\operatorname{mod}}
\nc{\End}{\operatorname{End}} \nc{\wh}[1]{\widehat{#1}} \nc{\Ext}{\operatorname{Ext}} \nc{\ch}{\text{ch}} \nc{\ev}{\operatorname{ev}}
 \makeatletter
\def\mydggeometry{\makeatletter\dg@YGRID=1\dg@XGRID=20\unitlength=0.003pt\makeatother}
\makeatother \theoremstyle{remark}

\numberwithin{equation}{section}
\let\bwdg\bigwedge
\def\bigwedge{{\textstyle\bwdg}}

\newcommand{\rnc}{\renewcommand}

\nc{\cal}{\mathcal} \nc{\goth}{\mathfrak} \rnc{\bold}{\mathbf}

\nc\bomega{{\mbox{\boldmath $\omega$}}} \nc\bpsi{{\mbox{\boldmath $\Psi$}}}
 \nc\balpha{{\mbox{\boldmath $\alpha$}}}
 \nc\bpi{{\mbox{\boldmath $\pi$}}}
 \nc\bvpi{{\mbox{\boldmath $\varpi$}}}

\newcommand\Lg{\mathfrak{g}}
\newcommand\Lh{\mathfrak{h}}
\newcommand\Lb{\mathfrak{b}}
\newcommand\Ln{\mathfrak{n}}

\DeclareMathOperator{\spa}{span}

\nc{\C}{\mathbb C }
\nc{\Z}{\mathbb Z }
\nc{\N}{\mathbb N }
\nc{\R}{\mathbb R }
\nc{\Q}{\mathbb Q }
\nc\blambda{{\mbox{\boldmath $\lambda$}}}
\nc\bmu{{\mbox{\boldmath $\mu$}}}
\nc\bsigma{{\mbox{\boldmath $\sigma$}}}

\newcommand{\gr}{\operatorname{gr}}

\DeclareMathOperator{\D}{D}

\DeclareMathOperator{\W}{W}

\begin{document}

\title[Representations of Lie superalgebras with fusion flags]{Representations of Lie superalgebras with fusion flags}

\author{Deniz Kus}
\address{Deniz Kus:\newline
Mathematisches Institut, Universit\"at zu K\"oln, Germany}
\email{dkus@math.uni-koeln.de}
\thanks{D.K was partially funded  under  the  Institutional  Strategy  of  the  University of Cologne within the German Excellence Initiative}

\subjclass[2010]{}
\begin{abstract}
We study the category of finite--dimensional representations for a basic classical Lie superalgebra $\Lg=\Lg_0\oplus \Lg_1$. For the ortho--symplectic Lie superalgebra $\Lg=\mathfrak{osp}(1,2n)$ we show that certain objects in that category admit a fusion flag, i.e. a sequence of graded $\Lg_0[t]$--modules such that the successive quotients are isomorphic to fusion products. Among these objects we find fusion products of finite--dimensional irreducible $\Lg$--modules, truncated Weyl modules and Demazure type modules. Moreover, we establish a presentation for these types of fusion products in terms of generators and relations of the enveloping algebra.
\end{abstract}
\maketitle \thispagestyle{empty}
\section{Introduction}
For a finite--dimensional simple Lie algebra $\mathfrak{t}$, the current algebra $\mathfrak{t}[t]$ associated to $\mathfrak{t}$ is the algebra of polynomial maps $\mathbb C \rightarrow \mathfrak{t}$ or equivalently, it is the complex vector space $(\mathfrak{t}\otimes \C[t])$ with Lie bracket the $\C[t]$--bilinear extension of the Lie bracket on $\mathfrak{t}$. The category of finite--dimensional representations of $\mathfrak{t}[t]$ has attracted a lot of attention over the past two decades. 
For example, motivated by the representation theory of quantum affine algebras the notion of local Weyl modules was introduced and studied in \cite{CP01}. Given any integrable highest weight representation of the affine Kac--Moody algebra $\widehat{\mathfrak{t}}$, one can define another family of finite--dimensional $\mathfrak{t}[t]$--modules. These modules are called Demazure modules and are parametrized by pairs $(\ell,\lambda)$, where $\ell$ is a positive integer called the level of the Demazure module and $\lambda$ is a dominant integral weight for $\mathfrak{t}$. It was proved in the simply--laced case (see \cite{CL06,CV13,FoL07}) and later for twisted current algebras (see \cite{CIK14, FK11, KV14}) that a local Weyl module is isomorphic to a level one Demazure module. In the non--simply laced case, a local Weyl module has a filtration by level one Demazure modules (see \cite{Na11,CSSW14}), i.e. there exists a flag where the successive quotients are isomorphic to Demazure modules. The multiplicity of a Demazure module in the local Weyl module has been studied in \cite{BCSV15,CSSW14}.\par
Another family of finite--dimensional representations for the current algebra was introduced in \cite{FL99}. Given finite--dimensional irreducible $\mathfrak{t}$--modules $V_{\mathfrak{t}}(\lambda_1),\dots, V_{\mathfrak{t}}(\lambda_k)$ and a tuple of pairwise distinct complex numbers $\mathbf z=(z_1,\dots,z_k)$ one can define a filtration on the tensor product $V_{\mathfrak{t}}(\lambda_1)^{z_1}\otimes \cdots \otimes V_{\mathfrak{t}}(\lambda_k)^{z_k}$ and build the associated graded space with respect to this filtration. This space is called the fusion product and is denoted by $V_{\mathfrak{t}}(\lambda_1)^{z_1}* \cdots *V_{\mathfrak{t}}(\lambda_k)^{*z_k}$, where $V_{\mathfrak{t}}(\lambda)^{z}$ is a non--graded $\mathfrak{t}[t]$--module. It was proved in \cite{CL06,FoL07,KL14} that certain (truncated) Weyl modules can be realized as fusion products of finite--dimensional irreducible $\mathfrak{t}$--modules.
We say that a representation $V$ has a fusion flag if there exists a sequence
$$\mathcal{F}(V)=\big(0\subset V_0\subset V_1\subset \cdots \subset V_k=V\big)$$
of graded $\mathfrak{t}[t]$--modules, such that the successive quotients are isomorphic to fusion products of finite--dimensional irreducible $\mathfrak{t}$--modules.\par
The present paper is motivated by the idea to study representations having a fusion flag. In particular, we show that the category of finite--dimensional representations of a Lie superalgebra contains many interesting objects with that property.\par
Let $\Lg=\Lg_0\oplus \Lg_1$ be a basic classical Lie superalgebra and let $\Lg[t]$ be the corresponding current superalgebra; recall that the even part $\Lg_0$ is a reductive Lie algebra and $\Lg_1$ is a semisimple $\Lg_0$--module. Starting with finite--dimensional irreducible $\mathfrak{g}$--modules $V(\lambda_1),\dots, V(\lambda_k)$ and a tuple of pairwise distinct complex numbers $\mathbf z=(z_1,\dots,z_k)$ one can define similarly as for finite--dimensional simple Lie algebras the fusion product $V(\lambda_1)^{z_1}* \cdots *V(\lambda_k)^{*z_k}$. For the ortho--symplectic Lie superalgebra $\Lg=\mathfrak{osp}(1,2n)$ we establish a presentation for certain fusion products in terms of generators and relations of the enveloping algebra. Using this presentation we show that fusion products for Lie superalgebras admit a fusion flag. Let us describe our results in more details for the Lie superalgebra $\mathfrak{osp}(1,2)$.\par We fix a partition $\mathbf{m}=(m_1\geq m_2\geq \cdots \geq m_k >m_{k+1}=0)$. For $0\leq \ell < k$ set
$$\varphi_{\ell}\big((m_1,\dots,m_k)\big):=(\varphi_{\ell}(m_1),\dots,\varphi_{\ell}(m_k))=(m_1,\dots,m_{p-1},m_p-1,m_{p+1},\dots,m_k),
$$ 
where $p=\min\big\{\ell+1 \leq q \leq k\mid m_q>m_{q+1}\big\}$. Our first results is the following; for the precise definition of the ingredients and a more general formulation see Section~\ref{section4}.
\begin{theom}\mbox{}
\begin{enumerate}
\item We have an isomorphism of $\mathbf{U}(\Lg[t])$--modules
$$V(m_1)*\cdots*V(m_k)\cong \W(|\mathbf{m}|)/\mathfrak{K}(\mathbf{m}),$$
where $\W(|\mathbf{m}|)$ is the local Weyl module and $\mathfrak{K}(\mathbf{m})$ is the submodule of $\W(|\mathbf{m}|)$ generated by the elements \eqref{7}.
\item The module $W(|\mathbf{m}|)/\mathfrak{K}(\mathbf{m})$ can be filtered by 
$$\bigoplus_{0\leq\ell\leq k}\ \bigoplus_{0\leq i_1 <\cdots <i_\ell\leq k-1} V_{\mathfrak{sl}_2}\big(\varphi_{i_1}\circ \cdots \circ \varphi_{i_\ell}(m_1)\big)*\cdots *V_{\mathfrak{sl}_2}\big(\varphi_{i_1}\circ \cdots \circ \varphi_{i_\ell}(m_k)\big)$$
\end{enumerate}
\end{theom}
The above theorem is proven more general for the Lie superalgebra $\mathfrak{osp}(1,2n)$ in Theorem~\ref{mainthm} and the proof uses the presentation of fusion products of certain finite--dimensional irreducible $\Lg_0$--modules. This fact is of independent interest and is stated in Section~\ref{section6}.
The connection to truncated Weyl modules and Demazure type modules is described by the following theorem. Again for the precise definition of the ingredients see Section~\ref{section5}.
\begin{theom}
Let $n\in \mathbb Z_{+}$ and write $n=kN+j$ for $0\leq j <N$. 
We have an isomorphism of $\mathbf U(\Lg[t])$--modules
\begin{enumerate}
\item $\W(n,N)\cong V(k)^{*(N-j)}*V(k+1)^{*j}$
\item $\W(n)\cong \D(1,n)$
\item $\D(N,n)\cong V(N)^{*k}*V(j).$
\end{enumerate}
\end{theom}
The paper is organized as follows. In Section~\ref{section2} we recall the basic properties of Lie superalgebras and collect the needed results from \cite{K78}. In Section~\ref{section3} we recall the notion of fusion products and show that the dimension of the fusion product increases along an appropriate order. In the next Section we give generators and relations for the fusion product of certain $\mathfrak{osp}(1,2n)$--modules and show that they admit a fusion flag. In Section~\ref{section5} we discuss the connection of truncated Weyl modules and Demazure type modules with fusion products. In the last Section we give generators and relations for certain fusion products associated to $\Lg_0$--modules.

%
\section{Lie superalgebras}\label{section2}
\subsection{}
We denote the set of complex numbers by $\mathbb{C}$ and, respectively, the set of integers, non--negative integers, and positive integers  by $\Z$, $\Z_+$, and $\mathbb{N}$. Unless otherwise stated, all the vector spaces considered in this paper are $\mathbb{C}$--vector spaces and $\otimes$ stands for $\otimes_\mathbb{C}$.
\subsection{}
Let $\Lg=\Lg_0\oplus \Lg_1$ be a basic classical Lie superalgebra and fix a Cartan subalgebra $\Lh_0\subset \Lg$, which is by definition a Cartan subalgebra of the even subalgebra $\Lg_0$. For $\alpha\in \Lh_0^{*}$ let
$$\Lg_{\alpha}=\big\{x\in \Lg\mid [h,x]=\alpha(h)x,\ \forall  h\in \Lh_0\big\}$$
the root space associated to $\alpha$ and $R=\big\{\alpha\in \Lh_0^{*}\mid \alpha\neq 0, \Lg_{\alpha}\neq 0\big\}$ be the root system of $\Lg.$
We define the even and odd roots to be
$$R_0=\big\{\alpha\in R \mid \Lg_{\alpha}\cap \Lg_0\neq 0\big\},\ R_1=\big\{\alpha\in R \mid \Lg_{\alpha}\cap \Lg_1\neq 0\big\}.$$
Similarly to semisimple Lie algebras, we have a root space decomposition 
\begin{equation}\label{roots}\Lg=\Lh_0\oplus \bigoplus_{\alpha\in R}\Lg_{\alpha},\end{equation}
where each root space $\Lg_{\alpha}$ in \eqref{roots} is one--dimensional.  We denote by $\W_0$ the Weyl group of $\Lg$, which is by definition the Weyl group of the reductive Lie algebra $\Lg_0$. There exists a non-degenerate even invariant supersymmetric bilinear form $(\cdot,\cdot)$ on $\Lg$, whose restriction to $\Lh_0\times \Lh_0$ is non--degenerate and $\W_0$--invariant. For a root $\alpha\in R$ we have $k\alpha\in R$ for $k\neq \pm 1$ if and only if $\alpha\in R_1$ and $(\alpha,\alpha)\neq 0$ (in this case $k=\pm 2)$.
\subsection{}
Let $E$ be the real vector space spanned by $R$ and $\succ$ a total ordering on $E$ compatible with the real vector space structure. We denote by $R^{+}=\{\alpha\in R \mid \alpha \succ 0\}$ and $R^{-}=\{\alpha\in R \mid \alpha \prec 0\}$ respectively the set of positive roots and negative roots respectively. We fix a subset $\Delta=\{\alpha_1,\dots,\alpha_r\}\subseteq R^+$ of simple roots, which by definition means that $\alpha\in \Delta$ cannot be written as a sum of two positive roots. We denote by $I=\{1,\dots,r\}$ the corresponding index set. Let $\rho_0$ (respectively $\rho_1$) be the half--sum of all the even (respectively odd) positive roots and set $\rho=\rho_0-\rho_1$.
We have a triangular decomposition 
$$\Lg=\Ln^-\oplus \Lh_0 \oplus \Ln^+,\mbox{ where } \Ln^{\pm}=\bigoplus_{\alpha\in R^{\pm}} \Lg_{\alpha}.$$
The subalgebra $\Lb=\Lh_0\oplus \Ln^+$ is solvable and is called the Borel subalgebra of $\Lg$ corresponding to the positive system $R^+$. We set $\Ln^{\pm}_{i}=\Ln^{\pm}\cap \Lg_i$ for $0\leq i \leq 1$.

\subsection{}
A postive root system is called distinguished if the corresponding system of simple roots contains exactly one odd root. \textit{From now on we fix a distinguished positive root system for $\Lg$ with Cartan matrix $A=(a_{i,j})_{i,j\in I}$ whose Dynkin diagram $S$ is given as in \cite[Table 1]{K78}. We denote by $s$ the unique node such that $\alpha_s$ is odd.}
Recall that the Cartan matrix $A$ of $\Lg$ satsifies the conditions
$$\begin{cases}a_{i,i}\in\{0,2\},\ \text{ for all $i=1,\dots,r$}\\ 
\text{if $a_{i,i}=0$, then $a_{i,i+k}=1$ where $k=\min\{1\leq j \leq n-i\mid a_{i,i+j}\neq 0\}$}.\end{cases}$$

Let $D=\mbox{diag}(d_i)_{i\in I}$ and $B=(b_{i,j})_{i,j\in I}$ be diagonal and symmetric matrices such that $A=DB$ (these matrices are explicitly given in \cite[Appendix]{IK01}). For the rest of this subsection we recall the notion of a Chevalley basis; for more details we refer to \cite{FG11,IK01}. For each positive root $\alpha\in R^+$ we fix a non--zero generator $x_{\alpha}^{\pm}$ of $\Lg_{\pm \alpha}$ and a linearly independent subset $\{h_i \mid i\in I\}$ of the Cartan subalgebra such that $\alpha_i(h_j)=a_{j,i}$. For a root $\alpha=\sum_{i=1}^rk_i\alpha_i\in R$ we define its coroot 
$$h_{\alpha}= d_{\alpha}\sum_{i=1}^r k_id_i^{-1}h_i\in \Lh_0,$$

where 
$$d_{\alpha}=\begin{cases}
\frac{2}{(\alpha,\alpha)},& \text{if $(\alpha,\alpha)\neq 0$}\\
d_s,& \text{if $(\alpha,\alpha)=0$.}
\end{cases}$$

We have
$$\alpha(h_{\alpha})=\begin{cases}2,& \text{ if $(\alpha,\alpha)\neq 0$}\\
0,& \text{ if $(\alpha,\alpha)=0$.}
\end{cases}
$$

The basis $\big\{x_{\alpha}^{\pm},\ h_i\mid i\in I,\ \alpha\in R^+\big\}$ of $\Lg$ is called a Chevalley basis if the following properties are satisfied:
\begin{equation}\label{0}\big\{h_1,\dots,h_r\big\} \mbox{ is a basis of }\Lh_0\end{equation}
\begin{equation}\label{1} [h_i,h_j]=0,\quad [h_i,x^{\pm}_{\alpha}]=\pm\alpha(h_i)x^{\pm}_{\alpha},\quad [x_{\alpha}^+,x_{\alpha}^-]=h_{\alpha},\quad \forall i,j\in I,\ \alpha\in R^+\end{equation}
\begin{equation}\label{2}[x_{\alpha},x_{\beta}]=c_{\alpha,\beta}x_{\alpha+\beta}, \quad \forall \alpha,\beta\in R \mbox{ with $\alpha\neq -\beta$ and $c_{\alpha,\beta}\in \mathbb{Z}\backslash \{0\}$}
\end{equation}
where $$x_{\alpha}:=\begin{cases}x_{\alpha}^{+},& \text{ if $\alpha\in R^+$}\\
 x_{\alpha}^{-},& \text{ if $\alpha\in R^-$}\\
0,& \text{ if $\alpha\notin R$}
\end{cases}
$$
The structure constants $c_{\alpha,\beta}$ satisfy further conditions; see for example \cite[Definition 3.3]{FG11}. For the rest of this paper, we fix such a Chevalley basis.

\subsection{}
We shall write the root system for the ortho--symplectic Lie superalgebra $\Lg=\mathfrak{osp}(2m+1,2n)$ ($m\geq 0, n\geq 1$). It consists of all matrices of the form

$$ 
\begin{pmatrix}\begin{array}{ccc|cc}
0 & -w^{t} & -v^{t} & x & x_1\\
v & A & B & Y& Y_1\\
w & C & -A^{t} & Z & Z_1\\
\hline
-x_1^{t} & -Z_1^{t} & -Y_1^{t} & D & G\\
x^{t} & Z^{t} & Y^{t} & F & -D^{t}
\end{array}
\end{pmatrix}
$$

where $A,B,C\in \mathbb{C}^{m\times m}$, $D,G,F\in \mathbb{C}^{n\times n}$, $Y,Y_1,Z,Z_1\in \mathbb{C}^{n\times m}$, $u,v\in \mathbb{C}^{m\times 1}$ and $x,x_1\in \mathbb{C}^{1\times n}$. The matrices $B,C$ are skew symmetric and the matrices $G,F$ are symmetric. We enumerate the rows as follows: the first row (resp. column) is indexed with 0 and the remaining rows (resp. columns) are indexed with the set $\{1,\dots,m,-1,\dots,-m,m+1,\dots,m+n,-(m+1),\dots,-(m+n)\}$. Let 
$$H_i=E_{i,i}-E_{-i,-i},\ 1\leq i \leq m, \quad H_j^{'}=E_{m+j,m+j}-E_{-(m+j),-(m+j)},\ 1\leq j\leq n,$$
then
$$\Lh_0=\spa\big\{H_1,\dots,H_m,H_1^{'},\dots,H_n^{'}\big\}$$
is a Cartan subalgebra of $\mathfrak{osp}(2m+1,2n)$. Let $\epsilon_1,\dots, \epsilon_m,\delta_1,\dots,\delta_n$ be elements in $\Lh_0^{*}$ defined by 
$$\epsilon_i(H_j)=\delta_{i,j},\ \epsilon_i(H_j^{'})=0,\ \delta_i(H_j^{'})=\delta_{i,j},\ \delta_{i}(H_j)=0.$$
The root system of $\mathfrak{osp}(2m+1,2n)$ is given by
$$R_0=\big\{\pm\epsilon_k\pm\epsilon_l,\ \pm\epsilon_q,\ \pm 2\delta_p,\ \pm \delta_i\pm \delta_j\big\},\quad R_1=\big\{\pm \delta_p,\ \pm\epsilon_q\pm \delta_p\big\},$$
where $1\leq k <l\leq m$, $1\leq i < j\leq n$, $1\leq p \leq n$ and $1\leq q \leq m$.
For the remaining part of this subsection we choose a distinguished positive root system
$$R^+=\Big\{\epsilon_k\pm\epsilon_l,\ \epsilon_q,\ 2\delta_p,\ \delta_i\pm \delta_{j}\Big\}\cup\Big\{\delta_p,\ \delta_p\pm\epsilon_q \Big\}$$
where the corresponding set of distinguished simple roots is given by
$$\begin{cases}\delta_1-\delta_{2},\dots,\delta_{n-1}-\delta_n,\ \delta_n-\epsilon_1,\ \epsilon_1-\epsilon_{2},\dots,\epsilon_{m-1}-\epsilon_{m},\ \epsilon_m & \text{ if $m>0$}\\
\delta_1-\delta_{2},\dots,\delta_{n-1}-\delta_n,\ \delta_n & \text{ else}.\end{cases}$$

Note that the only simple odd root is $\alpha_n$ and hence $R^+$ is distinguished. The elements $h_i, i\in I$ are given by
$$
\begin{cases}h_r=H_r^{'}-H_{r+1}^{'},\ 1\leq r <n,\ h_{n+s}=H_{s}-H_{s+1},\ 1\leq s<m,\ h_n=H_n^{'}+H_1,\ h_{n+m}=2 H_m& \text{ if $m>0$}\\
h_r=H_r^{'}-H_{r+1}^{'},\ 1\leq r <n,\ h_{n}=2 H_n^{'}& \text{ else.} \end{cases}
$$
For $m>0$ and each simple root we fix a generator of the root space (see \cite[pg. 21]{Mu12}): 
\begin{align}
x^+_{\alpha_r}&\propto E_{m+r,m+r+1}-E_{-(m+r+1),-(m+r)},\ 1\leq r <n \\
\notag x^-_{\alpha_r}&\propto E_{m+r+1,m+r}-E_{-(m+r),-(m+r+1)},\ 1\leq r <n\\
\notag x^+_{\alpha_{n+s}}&\propto E_{s,s+1}-E_{-(s+1),-s},\ x^-_{\alpha_{n+s}}=E_{s+1,s}-E_{-s,-(s+1)},\ 1\leq s<m \\
\notag x^+_{\alpha_{n}}&\propto E_{m+n,1}-E_{-1,-(m+n)},\ x^-_{\alpha_{n}}=E_{-(m+n),-1}+E_{1,m+n},\\
\notag x^+_{\alpha_{n+m}}&\propto E_{m,0}-E_{0,-m},\ x^-_{\alpha_{n+m}}=E_{-m,0}-E_{0,m}.
\end{align}
If $m=0$ we fix 
\begin{align}
x^+_{\alpha_r}&\propto E_{m+r,m+r+1}-E_{-(m+r+1),-(m+r)},\ 1\leq r <n \\
\notag x^-_{\alpha_r}&\propto E_{m+r+1,m+r}-E_{-(m+r),-(m+r+1)},\ 1\leq r <n\\
\notag x^+_{\alpha_{n}}&\propto E_{n,0}-E_{0,-n},\ x^-_{\alpha_{n}}=E_{-n,0}+E_{0,n}.
\end{align}
The set $\big\{x^{\pm}_{\alpha_i}, h_i\mid 1\leq i \leq n+m\big\}$ generates $\mathfrak{osp}(2m+1,2n)$ and satisfies the relations \eqref{1}.

\subsection{}
For a Lie superalgebra $\mathfrak{a}=\mathfrak{a}_0\oplus \mathfrak{a}_1$, we let $\mathbf{U}(\mathfrak{a})$ be the universal enveloping algebra of $\mathfrak{a}$. Let $B_0$ be a basis of $\mathfrak{a}_0$ and $B_1$ be a basis of $\mathfrak{a}_1$. If $\leq$ is a total order on $B=B_0\cup B_1$, then we obtain by the PBW theorem that the set of monomials
$$\Big\{g_1\cdots g_r\mid \mbox{ $g_j\in B,\ g_j\leq g_{j+1}$, and $g_j\neq g_{j+1}$ if $g_j\in \mathfrak{a}_1$}\Big\}$$
forms a basis of $\mathbf{U}(\mathfrak{a})$. The current superalgebra associated to $\mathfrak{a}$ is defined by $\mathfrak{a}[t]:=(\mathfrak{a}\otimes \mathbb{C}[t])$ where the even and odd part respectively is given by $\mathfrak{a}_0[t]$ and $\mathfrak{a}_1[t]$ respectively. The Lie superbracket is given in the obvious way. 

\subsection{}
For $\lambda\in \Lh_0^{*}$ we define a one--dimensional irreducible $\Lb$--module $\mathbb{C}_{\lambda}=\mathbb{C}v_\lambda$ by

$$\Ln^+v_{\lambda}=0,\ hv_{\lambda}=\lambda(h)v_{\lambda}, \ \forall h\in \Lh_0.$$
Consider the induced module $L(\lambda)=\mathbf{U}(\Lg)\otimes_{\mathbf{U}(\Lb)}\mathbb{C}_{\lambda}$ and let $J(\lambda)$ be the unique maximal submodule. We set $V(\lambda)=L(\lambda)/J(\lambda)$ and obtain that $V(\lambda)$ is an irreducible representation with $V=\mathbf{U}(\Ln^-)v_{\lambda}$ (for simplicity we denote the highest weight vector $1\otimes v_{\lambda}$ also by $v_{\lambda}$). Hence $V(\lambda)$ has a weight space decomposition
\begin{equation}\label{3}V(\lambda)=\bigoplus_{\mu\in \Lh_0^{*}}V_\mu,\quad V_{\mu}=\big\{v\in V \mid hv=\mu(h)v, \forall h\in \Lh_0\big\}\end{equation}
and $V_{\mu}\neq 0$ implies $\lambda-\mu$ is a $\mathbb{Z}_+$--linear combination of positive roots. The following proposition is proved in \cite[Proposition 2.2]{K78}.
\begin{prop}\mbox{}
\begin{enumerate}
\item We have $V(\lambda)\cong V(\mu)$ if and only if $\lambda=\mu$.
\item All weight spaces $V_{\mu}$ in \eqref{3} are finite--dimensional.
\item Any finite--dimensional irreducible $\Lg$--module is isomorphic to $V(\lambda)$ for some $\lambda\in \Lh_0^{*}$.
\end{enumerate}
\end{prop}
We are interested in finite--dimensional irreducible $\Lg$--modules. Let $P^+=\{\lambda\in \Lh_0^{*}\mid \dim V(\lambda)<\infty\}$.
\subsection{}
One of the neseccary conditions that $V(\lambda)$ is finite--dimensional is that $\lambda$ is a dominant integral weight for the Lie algebra $\Lg_0$. For the special linear Lie superalgebra $\mathfrak{sl}(m,n)$ this condition is also sufficient, i.e. $\lambda\in P^+$ if and only if $\lambda(h_i)\in \mathbb{Z}_+$ for all $i\neq s$. The characterizing properties for the remaining basic classical Lie superalgebras can be found in \cite[Proposition 2.3]{K78}. Recall that a representation $V(\lambda)$ is called typical if $(\lambda+\rho,\alpha)\neq 0$ for all $\alpha$ isotropic. \textit{From now on we assume that $\lambda\in P^+$ and $V(\lambda)$ is typical.} The following proposition stated in \cite[Theorem 1]{K78} gives generators and relations for typical finite--dimensional irreducible $\Lg$--modules.  
\begin{prop}\label{prop1}
We have an isomorphism of $\Lg$--modules
$$V(\lambda)\cong \mathbf{U}(\Lg)/M(\lambda),$$
where $M(\lambda)$ is the ideal generated by 
$$\begin{cases} \Ln^+,\ (h-\lambda(h))\mbox{ for all $h\in \Lh_0$},\ (x_{\alpha_i}^-)^{\lambda(h_i)+1}\mbox{ for } i\neq s,& \text{ if $\Lg$ is of type $A(m,n)$ or $C(n)$}\\
\Ln^+,\ (h-\lambda(h))\mbox{ for all $h\in \Lh_0$},\ (x_{\alpha_i}^-)^{\lambda(h_i)+1} \mbox{ for } i\neq s,\ (x^-_{\gamma})^{2\frac{(\lambda,\gamma)}{(\gamma,\gamma)}+1},& \text{ else}
\end{cases}
$$
and $\gamma=\sum_{i=s}^{r}c_i\alpha_i$ with labels $c_i$ as in \cite[Table 2]{K78}. Moreover,
$$\dim V(\lambda)= 2^{|R_1^+|}\prod_{\alpha\in R_0^+}\frac{(\lambda+\rho,\alpha)}{(\rho_0,\alpha)}.$$
\end{prop}
\begin{example}\mbox{}
\begin{enumerate}
\item
Consider the Lie superalgebra $\mathfrak{osp}(1,2)$. We have $\Lh_0=\spa\big\{h_1\big\}$, $R^+=\{2\delta_1\}\cup \{\delta_1\}$ and $\Delta=\{\delta_1\}$. The finite--dimensional irreducible representations are parametrized by the non--negative intergers. To be more precise, for $\lambda=m\delta_1$ we have $V(\lambda)$ is finite--dimensional if and only if $m\in \mathbb{Z}_+$ and 
$$V(\lambda)\cong \mathbf{U}(\Lg)/M(\lambda),$$
where $M(\lambda)$ is the ideal generated by $\Ln^+,\ (h-\lambda(h)),\mbox{ for all $h\in \Lh_0$},\ (x^-_{2\delta_1})^{m+1}$.
The dimension is given by $\dim V(m\delta_1)=2m+1.$
\item Consider the special linear Lie superalgebra $\mathfrak{sl}(2,1)$. We have $\Lh=\spa\big\{h_1,h_2\big\}$, $R^+=\{\epsilon_1-\epsilon_2\}\cup \{\epsilon_1-\delta_1,\epsilon_2-\delta_1\}$ and $\Delta=\{\epsilon_1-\epsilon_2,\epsilon_2-\delta_1\}$. Let $\{\omega_1,\omega_2\}$ be the dual basis of $\{h_1,h_2\}$ and $\lambda=m\omega_1+n\omega_2$. Then the irreducible representation $V(\lambda)$ is finite--dimensional if and only if $m\in\mathbb{Z}_+$ and typical if and only if $n+m\notin \{-1,0\}$. In this case we have
$$V(\lambda)\cong \mathbf{U}(\Lg)/M(\lambda),$$
where $M(\lambda)$ is the ideal generated by $\Ln^+,\ (h-\lambda(h)),\mbox{ for all $h\in \Lh_0$},\ (x^-_{\epsilon_1-\epsilon_2})^{m+1}$. 
\end{enumerate}
\end{example}
%

\section{Posets, tensor products and fusion products}\label{section3}
In this section we shall give the analogue result of \cite[Theorem 1 (i)]{CFS12} for the typical finite--dimensional irreducible representations for basic classical Lie superalgebras. The proof proceeds similarly and uses the dimension formula stated in Proposition~\ref{prop1}. We will use this fact as a motivation to study fusion products and prove certain surjective maps among them.
\subsection{}
We fix $\lambda\in P^+$ and let $P^+(\lambda,k)$ be the set of $k$--tuples $\blambda=(\lambda_1,\dots,\lambda_k)$, such that $\lambda_i\in P^+$, $1\leq i \leq k$ and $\sum_{j=1}^k \lambda_j=\lambda$. Let $\blambda=(\lambda_1,\dots,\lambda_k), \bmu=(\mu_1,\dots,\mu_k)\in P^+(\lambda,k)$. For a positive even root $\beta$ we define
$$r_{\beta,\ell}(\blambda)=\min\big\{(\lambda_{i_1}+\cdots+\lambda_{i_\ell})(h_{\beta})\mid 1\leq i_1<\cdots<i_\ell\leq k\big\}.$$
Note that definition makes sense since $\lambda(h_{\beta})\in \mathbb{Z}_+$ for all $\beta\in R_0^+$.
We say $\blambda \preceq \bmu$ if
\begin{equation}\label{4}r_{\beta,\ell}(\blambda)\leq r_{\beta,\ell}(\bmu) \mbox{ for all $\beta\in R_0^{+}$ and } 1\leq \ell\leq k.\end{equation}
We define an equivalence relation $\sim$ on $P^+(\lambda,k)$ by
$$\blambda\sim \bmu \Leftrightarrow \mbox{ $r_{\beta,\ell}(\blambda)=r_{\beta,\ell}(\bmu)$ for all $\beta\in R_0^{+}$ and $1\leq \ell\leq k$.}$$
The above partial order was introduced in \cite{CFS12}, where the authors made the following observation: for a finite--dimensional simple Lie algebra $\mathfrak{t}$ and a tuple $\blambda$ the dimension of the tensor product of the corresponding finite--dimensional irreducible $\mathfrak{t}$--modules increases along $\preceq$. We shall prove the analogue result for Lie superalgebras.
\subsection{} We will need the following lemma first.
\begin{lem}\label{5}
We have 
$2(\lambda+\rho)(h_\alpha)\in \mathbb{N}$ for all $\alpha\in R_0^+$ and $\lambda\in P^+$.
\proof
The proof is a case-by-case consideration. We know that a set of simple roots for $\Lg_0$ is given by $\{\alpha_i\mid i\neq s\}$ if $\Lg$ is of type $A(m,n)$ or $C(n)$ and $\{\alpha_i, \gamma \mid i\neq s\}$ otherwise, where $\gamma$ is as in Proposition~\ref{prop1}. Let $\alpha\in R_0^+$ and write $\alpha=\sum_{i\neq s}k_i\alpha_i+k\gamma$ (we assume that $k=0$ if the Lie superalgebra is of type $A(m,n)$ or $C(n)$). The corresponding coroot is given by 
$$h_{\alpha}=d_{\alpha}\sum_{i\neq s}d_i^{-1}k_ih_i+\frac{d_{\alpha}}{d_{\gamma}}kh_{\gamma}.$$
Since $d_{\alpha}d_i^{-1}k_i\in \mathbb{Z}_+$, $\lambda(h_i)\in \mathbb{Z}_+$ and $\rho(h_i)=1$ for all $i\neq s$, the claim follows for type $A(m,n)$ or $C(n)$. In the remaining types it will be sufficient to show $(\lambda+\rho)(h_\gamma)\in \mathbb{N}$, since an easy calculation shows that $\frac{d_{\alpha}}{d_{\gamma}}\in \mathbb{N}$ whenever $k\neq0$. We have
$(\lambda+\rho)(h_{\gamma})=\lambda(h_{\gamma})+d_{\gamma}(\rho,\gamma)$, where
$$\big((\rho,\gamma); (\gamma,\gamma)\big)=\begin{cases}\big(-2(m-1)-1;4\big),& \text{if $\Lg$ is of type $B(m,n)$}\\
\big(-2(m-1);4\big),& \text{if $\Lg$ is of type $D(m,n)$} \\
\big(9;-6\big),& \text{if $\Lg$ is of type $F(4)$}\\
\big(10;-8\big),& \text{if $\Lg$ is of type $G(3)$}\\
\big((1+\alpha);-2(1+\alpha)\big),& \text{if $\Lg$ is of type $D(2,1;\alpha)$.}
\end{cases}
$$
Since $\lambda\in P^+$ is typical we obtain that $\lambda(h_{\gamma})\geq b$, where $b$ is as in \cite[Table 2]{K78}. Thus in all cases we have that $(\lambda+\rho)(h_{\gamma})$ is a positive half integer. 
\endproof
\end{lem}

\begin{thm}
Let $\blambda\preceq \bmu$, then 
$$\prod_{i=1}^{k}\dim V(\lambda_i)\leq \prod_{i=1}^{k}\dim V(\mu_i)$$
and we have equality if and only if $\blambda=\bmu$ in $P^+(\lambda,k)/\sim$.
\proof
The statement in the theorem is equivalent to
$$2^{|R_1^+|}\prod_{\alpha\in R_0^+}\prod_{j=1}^{k}\frac{(\lambda_j+\rho,\alpha)}{(\rho_0,\alpha)}\leq 2^{|R_1^+|}\prod_{\alpha\in R_0^+}\prod_{j=1}^{k}\frac{(\mu_j+\rho,\alpha)}{(\rho_0,\alpha)}.$$
It suffices to show that
$$\prod_{j=1}^{k}2(\lambda_j+\rho)(h_\alpha)\leq \prod_{j=1}^{k}2(\mu_j+\rho)(h_\alpha), \mbox{ for each $\alpha\in R_0^+$.}$$
From Lemma~\ref{5} we know that $2(\lambda+\rho)(h_\alpha)\in \mathbb{N}$ for all $\alpha\in R_0^+$ and $\lambda\in P^+$. Hence we can order theses numbers 
$$0\leq 2(\lambda_{i_1}+\rho)(h_\alpha)\leq \cdots \leq 2(\lambda_{i_k}+\rho)(h_\alpha), \quad 0\leq 2(\mu_{j_1}+\rho)(h_\alpha)\leq \cdots \leq 2(\mu_{j_k}+\rho)(h_\alpha)$$
and $$\sum^{\ell}_{p=1}2(\lambda_{i_p}+\rho)(h_\alpha)\leq \sum^{\ell}_{p=1}2(\mu_{j_p}+\rho)(h_\alpha),\ \mbox{ for all $1\leq \ell \leq k$}.$$
The rest of the proof follows from an $\mathfrak{sl}_2$ argument as in \cite{CFS12} (see also the combinatorial result in the introduction).

\endproof
\end{thm}

The original motivation of \cite[Theorem 1 (i)]{CFS12} was to prove a stronger result, namely that $\blambda\preceq \bmu$ implies the existence of an injective homomorphism between the corresponding tensor products of finite--dimensional irreducible $\mathfrak{t}$--modules. Another way of showing the existence of these injective maps is to reformulate everything in the language of fusion products. We will recall the definition of fusion products in the next subsection; see \cite{FL99} for details.

\subsection{}
Let $V_1,\dots,V_k$ be finite--dimensional cyclic $\Lg$--modules with homogeneous cyclic vectors $v_1,\dots,v_k$. For a tuple of pairwise distinct complex numbers $\mathbf z= (z_1,\dots,z_k)\in \mathbb{C}^k$ consider the tensor product
\begin{equation}\label{tens}V_1^{z_1}\otimes \cdots \otimes V_k^{z_k},\end{equation}

where $V^{z_j}_j$ is a $\Lg[t]$--module whose action is given by 
$$(x\otimes f).w=f(z_j)x.w,\ x\in \Lg,\ f\in \mathbb{C}[t], \ w\in V_j.$$
The existence of the coproduct implies that the tensor product \eqref{tens} admits an $\mathbf{U}(\Lg[t])$ action, which is given by the following formula:
for $x\in \Lg$ homogeneous and $r\in \mathbb{Z}_+$
$$(x\otimes t^r)(v_1\otimes \cdots\otimes v_k)=\sum_{j=1}^{k} (-1)^{|x|(|v_1|+\cdots +|v_{j-1}|)}z_j^r(v_1\otimes\cdots \otimes x.v_j\otimes \cdots \otimes v_k),$$
where $|\cdot|$ denotes the parity.
The proof of the following lemma is the same of that in \cite[Proposition 1.6]{FL99} for finite--dimensional simple Lie algebras (see also \cite[Lemma 1.5]{FM15}).
\begin{lem}
For pairwise distinct complex numbers $z_1,\dots,z_k$, the tensor product $V^{z_1}_1\otimes \dots \otimes V^{z_k}_k$ is a cyclic $\mathbf{U}(\Lg[t])$--module with cyclic generator $v_{1}\otimes \cdots \otimes v_{k}$, i.e.
$$V^{z_1}_1\otimes \dots \otimes V^{z_k}_k\cong \mathbf{U}(\Lg[t])(v_{1}\otimes \cdots \otimes v_{k}).$$
\end{lem}
Recall that the ring $\mathbb{C}[t]$ is graded and hence we have an induced grading on $\mathbf{U}(\Lg[t])$. We denote by $\mathbf{U}(\Lg[t])[m]$ the homogeneous component of degree $m$ and note that it is a $\Lg$--module for all $m\in \mathbb{Z}_+$. We define an increasing filtration on the tensor product \eqref{tens}:
$$0\subset \mathbf{V}^0(\mathbf{z}):=\bigoplus_{m=0}^0\mathbf{U}(\Lg[t])[m](v_{1}\otimes \cdots \otimes v_{k}) \subset \mathbf{V}^1(\mathbf{z}):=\bigoplus_{m=0}^1\mathbf{U}(\Lg[t])[m](v_{1}\otimes \cdots \otimes v_{k}) \subset \cdots $$ 
of $\Lg$--modules.
The associated graded space 
$$\mbox{gr}\{V^j(\mathbf{z})\}_{j\geq 0}= \mathbf{V}^{0}(\mathbf{z})\oplus\bigoplus_{j\geq 1} \mathbf{V}^j(\mathbf{z}) /\mathbf{V}^{j-1}(\mathbf{z})$$ is called the fusion product and is denoted by $V^{z_1}_1* \dots * V^{z_k}_k$. We denote the image of $v_1\otimes\cdots \otimes v_k$ in the associated graded space by $v_1*\cdots * v_k$.
\begin{rem}
Note that the definition of the fusion product depends on the parameters $z_1,\dots,z_k$. It has been conjectured in \cite{FL99} that the fusion product of $\mathfrak{t}$--modules is, under suitable conditions on $V_j$ and $v_j$, independent of theses parameters. This conjecture has been proven in certain cases by various people (see for instance \cite{CL06,CSVW14,FF02,FoL07,KL14}). In this paper we will prove the independence for the fusion product of certain representations of the Lie superalgebra $\mathfrak{osp}(1,2n)$ and therefore we omit almost always the parameters in the notation.
\end{rem}

%
\section{Generators and Relations}\label{section4}
The aim of this section is to give generators and relations for the fusion product of certain finite--dimensional irreducible $\mathfrak{osp}(1,2n)$--modules. So let $\Lg=\mathfrak{osp}(1,2n)$ and recall that $\{\alpha_i=\delta_i-\delta_{i+1}, 1\leq i<n,\ \alpha_n=\delta_n\}$ is a set of distinguished simple roots. Recall that the even part $\Lg_0$ is isomorphic to the symplectic Lie algebra $\mathfrak{sp}_{2n}$ and the odd part is a $2n$--dimensional irreducible $\Lg_0$--module. For example, if $n=1$ the even part in spanned by $\{x_{2\delta_1}^+,x_{2\delta_1}^-,h_{2\delta_1}\}$ and is isomorphic to $\mathfrak{sl}_2$.

\subsection{} The definition of local Weyl modules was originally given in \cite{CP01}. For $k\in\mathbb{Z}_+$, the local Weyl module $\W_{\Lg_0}(k\delta_1)$ is the $\Lg_0[t]$--module generated by a non--zero element $w_k$ with defining relations:
$$\Ln^+_0[t]w_{k}=0,\ (h\otimes t^r)w_{k}=k\delta_1(h)\delta_{r,0}w_{k},\  \forall r\in \mathbb{Z}_+,\ h\in \Lh_0, \ (x_{\alpha}^{-}\otimes 1)^{k\delta_1(h_{\alpha})+1}w_{k}=0,\  \forall \alpha\in R_0^+.$$
The dimension of the local Weyl module $\W_{\Lg_0}(k\delta_1)$ for $n=1$ is computed in \cite{CP01} and equals $2^k$.
The study of local Weyl modules for Lie superalgebras has recently started; see \cite{CLS15,FM15}.
For $k\in \mathbb{Z}_+$, the local Weyl module $\W(k\delta_1)$ is the $\Lg[t]$--module generated by a non--zero element $w_k$ with defining relations:
$$\Ln^+[t]w_{k}=0,\ (h\otimes t^r)w_{k}=k\delta_1(h)\delta_{r,0}w_{k},\  \forall r\in \mathbb{Z}_+,\ h\in \Lh_0, \ (x_{\alpha}^{-}\otimes 1)^{k\delta_1(h_{\alpha})+1}w_{k}=0,\  \forall \alpha\in R_0^+.$$
The dimension of $\W(k\delta_1)$ for $n=1$ is computed in \cite{FM15} and equals $3^k$.

\subsection{}
For $s,r\in \Z_+$, we set
\begin{equation}\label{maxs}\mathbf S(r,s)=\Big\{(b_k)_{k\ge 0}\mid b_k\in\mathbb Z_+, \ \  \sum_{k\ge 0} b_k=r,\ \ \sum_{k\ge 0} kb_k=s\Big\}. \end{equation}
and for $p\in \Z_+$ let $_p\mathbf S(r,s)$ be the subset of $\mathbf S(r,s)$ consisting of all elements $(b_k)_{k\ge 0}$ satisfying $b_k=0$ for $k< p$. For an even positive root $\alpha$, we define elements of $\mathbf{U}(\Lg_0[t])$
\begin{equation}\label{xlm}\mathbf{x}^-_{\alpha}(r,s)=\sum_{(b_k)_{k\ge 0}\in\mathbf S(r,s)}(x_{\alpha}^-\otimes 1)^{(b_0)}(x_{\alpha}^-\otimes t)^{(b_1)}\cdots (x_{\alpha}^-\otimes t^s)^{(b_s)},\end{equation}
where $x^{(b)}:=\frac{1}{b!}x^b$ for any non--negative integer $b$ and $x\in \Lg[t]$. We understand that $\mathbf{x}^-_{\alpha}(r,s)$ is zero if $\mathbf S(r,s)$ is empty. The elements $_p\mathbf{x}^-_{\alpha}(r,s)$ are defined in the obvious way. 
\subsection{}
For the rest of this paper we fix a partition $\mathbf{m}=(m_1\geq m_2\geq \cdots \geq m_k>m_{k+1}=0)$ and let $V_{\Lg_0}(k\delta_1)$ be the $\binom{k+2n-1}{2n-1}$--dimensional irreducible $\Lg_0$--module. The following theorem is proved in Section~\ref{section6} and gives generators and relations for the fusion product of certain finite--dimensional irreducible $\Lg_0$--modules. The case $n=1$ is proven in \cite{CV13} and can be also deduced from \cite{FF02}.
\begin{thm}\label{appe}
We have an isomorphism of $\mathbf{U}(\Lg_0[t])$--modules
\begin{equation}\label{cnfall}V_{\Lg_0}(m_1\delta_1)*\cdots * V_{\Lg_0}(m_k\delta_1)\cong \W_{\Lg_0}(|\mathbf{m}|\delta_1)/\mathfrak{I}(\mathbf{m}),\end{equation}
where $\mathfrak{I}(\mathbf{m})$ is the submodule generated by the elements
$$\mathbf{x}^-_{\alpha}(r,s),\ \mbox{ for all }\alpha\in R^+_0,\ r,s\in \mathbb{N}\mbox{ with }\ s+r\ge 1+ rp+\sum_{j\ge p+1}m_j\delta_1(h_{\alpha}),\ \mbox{ for some $p\in \mathbb{Z}_+$}.$$
\end{thm}

We shall prove in the main theorem of this section that the similar relations determine the fusion product of finite--dimensional irreducible $\Lg$--modules. The following corollary follows from the above theorem and \cite[Section2]{CV13}, where three different presentations of the right hand side of \eqref{cnfall} are given; see also \cite{KV14} for the twisted analogues. 
\begin{cor}\label{appe2}
We have an isomorphism of $\mathbf{U}(\Lg_0[t])$--modules
$$\W_{\Lg_0}(|\mathbf{m}|\delta_1)/\mathfrak{I}(\mathbf{m})\cong \W_{\Lg_0}(|\mathbf{m}|\delta_1)/\mathfrak{N}(\mathbf{m}),$$
where $\mathfrak{N}(\mathbf{m})$ is the submodule generated by the elements
$$_p\mathbf{x}^-_{\alpha}(r,s),\ \mbox{ for all }\alpha\in R^+_0,\ r,s\in \mathbb{N},\ p\in \Z_+,\mbox{ with }\ s+r\ge 1+ rp+\sum_{j\ge p+1}m_j\delta_1(h_{\alpha}).$$
\end{cor}
\subsection{}\label{ord}
We define a total order on the set of multi--exponents. For $0\leq i_1<\cdots <i_\ell\leq k-1$ and $0\leq j_1<\cdots <j_s\leq k-1$ we say 
$$(j_{s},\dots,j_1) \prec (i_{\ell},\dots,i_1)$$
if $s< \ell$ or $s=\ell$ and $(i_{\ell},\dots,i_1)<(j_{\ell},\dots,j_1)$ with respect to the lexicographical order. We consider the induced total order on the set of monomials, i.e.
$$(x_{\delta_1}^{-}\otimes t^{j_1})\cdots (x_{\delta_1}^{-}\otimes t^{j_s}) \prec (x_{\delta_1}^{-}\otimes t^{i_1})\cdots (x_{\delta_1}^{-}\otimes t^{i_\ell}):\Leftrightarrow (j_{s},\dots,j_1) \prec (i_{\ell},\dots,i_1).$$ For $\ell\in \Z_+$, we set 
$$\text{G}_{i_1,\dots,i_{\ell}}=\sum_{(j_{s},\dots,j_1) \prec (i_{\ell},\dots,i_1)}\mathbf{U}(\Lg_0[t])(x_{\delta_1}^{-}\otimes t^{j_1})\cdots (x_{\delta_1}^{-}\otimes t^{j_s})
$$
if $l>0$ and $\text{G}_{i_1,\dots,i_{\ell}}=\text{G}_{\emptyset}=0$ otherwise.

\begin{lem}\label{brauchb2}
Let $V$ be a $\mathbf{U}(\Lg[t])$--module and $v\in V$ such that 
$$\Ln^+[t]v=\big(\Lh\otimes t\C[t]\big)v=0.$$
Then for all $0\leq i_1<\cdots <i_\ell\leq k-1$ we have
$$\Big(\Ln^+[t]\oplus \big(\Lh\otimes t\C[t]\big)\Big)(x_{\delta_1}^{-}\otimes t^{i_1})\cdots (x_{\delta_1}^{-}\otimes t^{i_\ell})v\in \text{G}_{i_1,\dots,i_{\ell}}v.$$
\proof
For $s\in \N$ we have
\begin{equation}\label{rt}(h\otimes t^s)(x_{\delta_1}^{-}\otimes t^{i_1})\cdots (x_{\delta_1}^{-}\otimes t^{i_\ell})v=\sum^{\ell}_{j=1} (x_{\delta_1}^{-}\otimes t^{i_1})\cdots (x_{\delta_1}^{-}\otimes t^{i_j+s})\cdots (x_{\delta_1}^{-}\otimes t^{i_\ell})v.\end{equation}
Since $[x_{\delta_1}^-,x_{\delta_1}^{-}]\in \Ln^-_0$ and $[\Ln_0^-,x_{\delta_1}^{-}]=0$ we can reorder each monomial in \eqref{rt} by the degree of $t$ and obtain
\begin{align*}(x_{\delta_1}^{-}\otimes t^{i_1})\cdots (x_{\delta_1}^{-}\otimes t^{i_{j}+s})\cdots (x_{\delta_1}^{-}\otimes t^{i_\ell})v=(x_{\delta_1}^{-}\otimes t^{j_1})\cdots (x_{\delta_1}^{-}\otimes t^{j_\ell})v+ \mbox{ some element in } \text{G}_{i_1,\dots,i_\ell}v.\end{align*}
We have $\{i_1,\dots,i_{j}+s,\dots, i_{\ell}\}=\{j_1,\dots,j_{\ell}\}$ and $j_1<\cdots <j_\ell$, which implies $$(x_{\delta_1}^{-}\otimes t^{j_1})\cdots (x_{\delta_1}^{-}\otimes t^{j_\ell})\in \text{G}_{i_1,\dots,i_\ell}.$$
The proof for the remaining elements works similar.
\endproof
\end{lem}
\begin{rem}
The above lemma and an easy induction argument implies 
$$\text{G}_{i_1,\dots,i_{\ell}}v=\sum_{(j_{s},\dots,j_1) \prec (i_{\ell},\dots,i_1)}\mathbf{U}(\Ln^-_0[t])(x_{\delta_1}^{-}\otimes t^{j_1})\cdots (x_{\delta_1}^{-}\otimes t^{j_s})v.$$
\end{rem}
\subsection{}
For $0\leq \ell < k$ set
$$\varphi_{\ell}\big((m_1,\dots,m_k)\big):=(\varphi_{\ell}(m_1),\dots,\varphi_{\ell}(m_k))=(m_1,\dots,m_{p-1},m_p-1,m_{p+1},\dots,m_k),
$$ 
where $p=\min\big\{\ell+1 \leq q \leq k\mid m_q>m_{q+1}\big\}$. 
We now prove,
\begin{prop}\label{brauchb}We fix $0\leq i_0< i_1<\cdots <i_\ell\leq k-1$ and $\alpha\in R_0^+$. Let $V$ be a $\mathbf{U}(\Lg[t])$--module and $v\in V$ such that 
\begin{equation}\label{rtt}\Big(\Ln^+[t]\oplus \big(\Lh\otimes t\C[t]\big)\oplus \sum^n_{j=2}(x_{\delta_j}^{-}\otimes \C[t])\Big)v=0\end{equation}
and
\begin{equation}\label{hilf0}_p\mathbf{x}_{\alpha}^{-}(r,s)\prod_{j=1}^{\ell}(x_{\delta_1}^{-}\otimes t^{i_j})v\in \text{G}_{i_1,\dots,i_{\ell}}v,\  \forall \ r,s\in \mathbb{N}, p\in \Z_+ : r+s\geq 1+rp+\sum_{j\geq p+1}m_j\delta_1(h_{\alpha}).\end{equation}
Then we get
\begin{equation}\label{hilf}_p\mathbf{x}_{\alpha}^{-}(r,s)\prod_{j=0}^{\ell}(x_{\delta_1}^{-}\otimes t^{i_j})v\in \text{G}_{i_0,\dots,i_{\ell}}v,\ \forall \ r,s\in \mathbb{N},\ p\in \Z_+: r+s\geq 1+rp+\sum_{j\geq p+1}\varphi_{i_0}(m_j)\delta_1(h_{\alpha}).\end{equation}
\proof
Let $r,s\in \mathbb{N}$, $p\in \Z_+$ such that $$r+s\geq 1+rp+\sum_{j\geq p+1}\varphi_{\ell}(m_j)\delta_1(h_{\alpha}).$$ Since $[x_{\alpha}^-,x_{\delta_1}^-]=0$ and 
$$(x_{\delta_1}^{-}\otimes t^{i_0})\text{G}_{i_1,\dots,i_{\ell}}v\subset \text{G}_{i_0,i_1,\dots,i_{\ell}}v,$$
the proposition is immediate for $p\geq k$ or $\delta_1(h_{\alpha})=0$. So we can suppose that $p<k$ and $\delta_1(h_{\alpha})=1$. This means that $\alpha$ is of the form $\alpha=\delta_1\pm \delta_j$ for some $2\leq j \leq n$ or $\alpha=2\delta_1$. Note that there exists a root $\beta$ such that the corresponding root vector $x_{\beta}$ satisfies $[x_{\beta},x_{\alpha}^-]=x_{\delta_1}^-$. To be more precise we have $x_\beta\propto x^{\pm}_{\delta_j}$ if $\alpha=\delta_1\pm\delta_j$ and $x_\beta\propto x^+_{\delta_1}$ otherwise. The proof proceeds in two steps, where the first step considers the case $p\leq i_0$. In this case we get 
$$(r+1)+(s+p)\geq 1+(r+1)p+\sum_{j\geq p+1}\varphi_{i_0}(m_j)+1=1+(r+1)p+\sum_{j\geq p+1}m_j$$ and hence by \eqref{hilf0} 
$$_p\mathbf{x}_{\alpha}^{-}(r+1,s+p)\prod_{j=1}^{\ell}(x_{\delta_1}^{-}\otimes t^{i_j})v\in \text{G}_{i_1,\dots,i_{\ell}}v.$$ We obtain 
\begin{align*}&\big(x_{\beta}\otimes t^{i_0-p}\big)\Big(_p\mathbf{x}_{\alpha}^{-}(r+1,s+p)\prod_{j=1}^{\ell}(x_{\delta_1}^{-}\otimes t^{i_j})\Big)v=&\\&
=\sum_{\substack{(b_k)_{k\geq p}\in _p\mathbf{S}(r+1,s+p)\\ b_p\neq 0}}(x_{\alpha}^{-}\otimes t^p)^{(b_p-1)}\cdots (x_{\alpha}^{-}\otimes t^{s})^{(b_{s})}\prod_{j=0}^{\ell}(x_{\delta_1}^{-}\otimes t^{i_j})v + \mbox{ some element in $\text{G}_{i_0,i_1,\dots,i_{\ell}}v$}
&\\&
=_p\mathbf{x}_{\alpha}^{-}(r,s)\prod_{j=0}^{\ell}(x_{\delta_1}^{-}\otimes t^{i_j})v + \mbox{ some element in $\text{G}_{i_0,i_1,\dots,i_{\ell}}v$}\in \big(x_{\beta}\otimes t^{i_0-p}\big)\text{G}_{i_1,\dots,i_{\ell}}v\subset \text{G}_{i_0,i_1,\dots,i_{\ell}}v.\end{align*}
This proves \eqref{hilf} for all $r,s\in \mathbb{N}$, $p\in \Z_+$ such that 
$$p\leq i_0 \mbox{ and } r+s\geq 1+rp+\sum_{j\geq p+1}\varphi_{\ell}(m_j)\delta_1(h_{\alpha}).$$ Assume now that $p>i_0$. In this case there is only something to prove if $m_{i_0+1}=\cdots=m_{p+1}$, because otherwise we have $\sum_{j\geq p+1}\varphi_{\ell}(m_j)=\sum_{j\geq p+1}m_j$ and \eqref{hilf} follows directly from \eqref{hilf0}. Here we also consider two cases, starting with $r\geq m_{p+1}$. We have
$$r+s\geq 1+rp+\sum_{j\geq p+1}\varphi_{i_0}(m_j)=rp+\sum_{j\geq p+1}m_j\geq r i_0+\sum_{j\geq i_0+1}m_j=1+ri_0+\sum_{j\geq i_0+1}\varphi_{i_0}(m_j)$$
and hence we can deduce from the first step
\begin{equation}\label{hilf2}_{i_0}\mathbf{x}_{\alpha}^{-}(r,s)\prod_{j=0}^\ell(x_{\delta_1}^{-}\otimes t^{i_j})v\in \text{G}_{i_0,i_1,\dots,i_{\ell}}v.\end{equation}
Our aim is to show by induction on $r$ that
\begin{equation}\label{hilf3}_{i_0+1}\mathbf{x}_{\alpha}^{-}(r,s)\prod_{j=0}^\ell(x_{\delta_1}^{-}\otimes t^{i_j})v\in \text{G}_{i_0,i_1,\dots,i_{\ell}}v.\end{equation}
If $r=1$, we have $_{i_0+1}\mathbf{x}_{\alpha}^{-}(1,s)=(x^-_{\alpha}\otimes t^s)=_{i_0}\mathbf{x}_{\alpha}^{-}(1,s)$ and the claim follows. Assume now that we have proved the statement for all $\tilde{r}<r$ and $\tilde{s}\leq s$ such that $\tilde{r}+\tilde{s}\geq \tilde{r}p+\sum_{j\geq p+1}m_j$.
Note that
\begin{align*}_{i_0}\mathbf{x}_{\alpha}^{-}(r,s)&\prod_{j=0}^\ell(x_{\delta_1}^{-}\otimes t^{i_j})v&\\&=_{i_0+1}\mathbf{x}_{\alpha}^{-}(r,s)\prod_{j=0}^\ell(x_{\delta_1}^{-}\otimes t^{i_j})v+\sum^{r}_{i=1}(x_{\alpha}^{-}\otimes t^{i_0})^{(i)}\Big(_{i_0+1}\mathbf{x}_{\alpha}^{-}(r-i,s-i_0 i)\Big)\prod_{j=0}^\ell(x_{\delta_1}^{-}\otimes t^{i_j})v,
\end{align*}
where we understand $_{i_0+1}\mathbf{x}_{\alpha}^{-}(r-i,s-i_0 i)=0$ if $_{i_0+1}\mathbf{S}(r-i,s-i_0 i)=\emptyset$. Since
$$(r-i)+(s-i_0 i)\geq rp+\sum_{j\geq p+1}m_j-i(i_0+1)\geq (r-i)p+\sum_{j\geq p+1}m_j$$
we deduce \eqref{hilf3} from the induction hypothesis and \eqref{hilf2}.
If $p=i_0+1$ we are done; otherwise we repeat the above procedure until we get \eqref{hilf}. Now it remains to consider the case $r<m_{p+1}$.
In that case we get
\begin{align*}_{p}\mathbf{x}_{\alpha}^{-}(r,s)\prod_{j=0}^\ell(x_{\delta_1}^{-}\otimes t^{i_j})v=\sum^{r}_{i=0}(x_{\alpha}^{-}\otimes t^{i_0})^{(i)}\Big(_{p+1}\mathbf{x}_{\alpha}^{-}(r-i,s-i_0 i)\Big)\prod_{j=0}^\ell(x_{\delta_1}^{-}\otimes t^{i_j})v,
\end{align*}
and $$(r-i)+(s-i_0 i)\geq rp+\sum_{j\geq p+1}m_j-i(i_0+1)\geq 1+(r-i)(p+1)+\sum_{j\geq p+2}m_j.$$ Hence \eqref{hilf0} implies $$_{p+1}\mathbf{x}_{\alpha}^{-}(r-i,s-i_0 i) \prod_{j=1}^\ell(x_{\delta_1}^{-}\otimes t^{i_j})v\in \text{G}_{i_1,\dots,i_{\ell}}v$$ and the proposition is proven.
\endproof
\end{prop}

\subsection{}
We denote by $\mathfrak{K}(\mathbf{m})$ the submodule of $\W(|\mathbf{m}|\delta_1)$ generated by the following elements
\begin{equation}\label{7}_p\mathbf{x}^-_{\alpha}(r,s),\ \mbox{ for all }\alpha\in R^+_0,\ s,r\in\mathbb{N},\ p\in \Z_+ \mbox{ with} \ s+r\ge 1+ rp+\sum_{j\ge p+1}m_j\delta_1(h_{\alpha}) \end{equation}
The proof of the next propostion proceeds exactly in the same way as for $\mathfrak{sl}_2$ fusion products; for details see \cite{CV13}.
\begin{prop}\label{surj}
We have a surjective map of $\mathbf{U}(\Lg[t])$--modules
 $$\W(|\mathbf{m}|\delta_1)/\mathfrak{K}(\mathbf{m})\twoheadrightarrow V(m_1\delta_1)*\cdots*V(m_k\delta_1).$$
\end{prop}
In particular, the proposition gives a lower bound for the dimension, namely 
\begin{equation}\label{10}\dim \W(|\mathbf{m}|\delta_1)/\mathfrak{K}(\mathbf{m})\geq \prod_{j=1}^k\dim V(m_j\delta_1).\end{equation} 
For the Lie superalgebra $\mathfrak{osp}(1,2)$ this means
$$\dim \W(|\mathbf{m}|\delta_1)/\mathfrak{K}(\mathbf{m})\geq \prod_{j=1}^k(2m_j+1).$$
\subsection{}
The main theorem of this section is the following.
\begin{thm}\label{mainthm}\mbox{}
\begin{enumerate}
\item We have an isomorphism of $\mathbf{U}(\Lg[t])$--modules
\begin{equation}\label{8}V(m_1\delta_1)*\cdots*V(m_k\delta_1)\cong \W(|\mathbf{m}|\delta_1)/\mathfrak{K}(\mathbf{m}).\end{equation}
\item The module $W(|\mathbf{m}|\delta_1)/\mathfrak{K}(\mathbf{m})$ can be filtered by 
\begin{equation}\label{9}\bigoplus_{0\leq\ell\leq k}\ \bigoplus_{0\leq i_1 <\cdots <i_\ell\leq k-1} V_{\Lg_0}\big(\varphi_{i_1}\circ \cdots \circ \varphi_{i_\ell}(m_1)\delta_1\big)*\cdots *V_{\Lg_0}\big(\varphi_{i_1}\circ \cdots \circ \varphi_{i_\ell}(m_k)\delta_1\big),\end{equation}
where the cyclic vectors are the images of the vectors in $$\Big\{(x_{\delta_1}^-\otimes t^{i_1})\cdots (x_{\delta_1}^-\otimes t^{i_\ell})v_\mathbf{m}\mid 0\leq \ell \leq k,\ 0\leq i_1 <\cdots <i_k\leq k-1\Big\}.$$
\end{enumerate}
\proof

We denote by $V(\mathbf{m})$ the right hand side of \eqref{8} and let $v_\mathbf{m}:= w_{|\mathbf{m}|} \mod \mathfrak{K}(\mathbf{m})$ be the cyclic generator. From the defining relations of the local Weyl module we know
$$(x_{2\delta_n}^-\otimes 1)v_\mathbf{m}=0 \Rightarrow (x^+_{\delta_n}\otimes 1)(x_{2\delta_n}^-\otimes 1)v_\mathbf{m}=\pm(x_{\delta_n}^-\otimes 1)v_\mathbf{m}=0.$$
Since $(x_{\delta_j-\delta_{j+1}}^-\otimes 1)v_\mathbf{m}=0$ for $2\leq j\leq n-1$ we get
$$(x_{\delta_j-\delta_{j+1}}^-\otimes 1)\cdots(x_{\delta_{n-1}-\delta_{n}}^-\otimes 1)(x_{\delta_n}^-\otimes 1)v_\mathbf{m}=(x^-_{\delta_j}\otimes 1)v_\mathbf{m}=0,\ \mbox{ for $2\leq j \leq n-1$}.$$
Together with the above calculations and $\big(\Lh\otimes t\mathbb{C}[t]\big)v_\mathbf{m}=0$ we deduce 
\begin{equation}\label{hhh}\big(x^-_{\delta_j}\otimes \C[t]\big)v_\mathbf{m}=0 \mbox{ for $2\leq j \leq n$}.\end{equation}
It follows that the cyclic generator $v_\mathbf{m}$ satisfies the condition \eqref{rtt}.
For $r=1$ and $s=k$ we have 
$$\mathbf{x}^-_{\alpha}(r,s)=(x_{\alpha}^{-}\otimes t^k)\in \mathfrak{K}(\mathbf{m})\mbox{ and thus }(x_{\delta_1}^+\otimes 1)(x_{2\delta_1}^{-}\otimes t^k)=(x_{2\delta_1}^{-}\otimes t^k)(x_{\delta_1}^+\otimes 1)\pm (x_{\delta_1}^-\otimes t^k)\in \mathfrak{K}(\mathbf{m}).$$ It follows that $\big(x_{\delta_1}^-\otimes t^k\C[t]\big)v_\mathbf{m}=0$ and hence by the PBW theorem we can write
$$V(\mathbf{m})=\sum_{\ell=0}^k\ \sum_{0\leq i_1<\cdots <i_\ell\leq k-1}\mathbf{U}(\Lg_0[t])(x_{\delta_1}^{-}\otimes t^{i_1})\cdots (x_{\delta_1}^{-}\otimes t^{i_\ell})v_{\mathbf{m}}.$$
Our aim is to show that $V(\mathbf{m})$ can be filtered by \eqref{9}. We order the finite set 
$$\Big\{(x_{\delta_1}^{-}\otimes t^{i_1})\cdots (x_{\delta_1}^{-}\otimes t^{i_\ell})\mid 0\leq \ell \leq k,\ 0\leq i_1<\cdots <i_\ell\leq k-1\Big\}=\Big\{X_1\prec X_2\prec \cdots \prec X_N\Big\}$$
with respect to the order defined in Section~\ref{ord}. We build the associated graded space with respect to the increasing filtration of $\Lg_0[t]$--modules
$$V(\mathbf{m})_1\subseteq V(\mathbf{m})_2\subseteq \cdots \subseteq V(\mathbf{m})_N=V(\mathbf{m}),$$
where $$V(\mathbf{m})_{i}=\sum_{j=1}^i\mathbf{U}(\Lg_0[t])X_jv_{\mathbf{m}}.$$
In particular we will show that we have a surjective map 
\begin{equation}\label{12}V_{\Lg_0}\big(\varphi_{i_1}\circ \cdots \circ \varphi_{i_\ell}(m_1)\delta_1\big)*\cdots *V_{\Lg_0}\big(\varphi_{i_1}\circ \cdots \circ \varphi_{i_\ell}(m_k)\delta_1\big)\twoheadrightarrow
V(\mathbf{m})_i/V(\mathbf{m})_{i-1},\end{equation}
where $X_i=(x_{\delta_1}^{-}\otimes t^{i_1})\cdots (x_{\delta_1}^{-}\otimes t^{i_\ell})$. We denote the image of $X_iv_{\mathbf{m}}$ in the associated graded space by $X_i\gr(v_{\mathbf{m}})$. For this it is enough to verify that $X_i\gr(v_{\mathbf{m}})$ satisfies the defining relations of the fusion product stated in Corollary~\ref{appe2}. From Lemma~\ref{brauchb2} we get that $X_i\gr(v_{\mathbf{m}})$ satisfies the defining relations of the local Weyl module $\W(|\mathbf{m}|\delta_1)$. Hence it remains so show that 
\begin{equation}\label{11}_p\mathbf{x}^-_{\alpha}(r,s)X_i\gr(v_{\mathbf{m}})=0,\ \forall s,r\in\mathbb{N},\ p\in \Z_+ : s+r\ge 1+ rp+\sum_{j\ge p+1}\varphi_{i_1}\circ \cdots \circ \varphi_{i_\ell}(m_j)\delta_1(h_{\alpha}),\end{equation}
which we will show by induction on $\ell$. If $\ell$ is zero, this follows from the definition of the submodule $\mathfrak{K}(\mathbf{m})$. So suppose that we have proved \eqref{11} for all $\widetilde{\ell}<\ell$. In particular, from the induction hypothesis 
$$_p\mathbf{x}^-_{\alpha}(r,s)\prod_{j=2}^{\ell}(x_{\delta_1}^{-}\otimes t^{i_j})\gr(v_{\mathbf{m}})=0,\ \forall s,r\in\mathbb{N},\ p\in \Z_+ : s+r\ge 1+ rp+\sum_{j\ge p+1}\varphi_{i_2}\circ \cdots \circ \varphi_{i_\ell}(m_j)\delta_1(h_{\alpha}).$$ 
By using Proposition~\ref{brauchb} we have
$$_p\mathbf{x}^-_{\alpha}(r,s)\prod_{j=1}^{\ell}(x_{\delta_1}^{-}\otimes t^{i_j})v_{\mathbf{m}}\in \text{G}_{i_1,\dots,i_{\ell}}$$
and hence the desired property \eqref{11}. An easy calculation shows
and $$\sum_{0\leq \ell \leq k-1}\ \sum_{0\leq i_1 <\cdots <i_k\leq k-1} \prod_{j=1}^k \dim V_{\Lg_0}\big(\varphi_{i_1}\circ \cdots \circ \varphi_{i_\ell}(m_j)\delta_1\big) =\prod_{j=1}^k\dim V(m_j\delta_1),$$
which implies together with \eqref{10} that the map in \eqref{12} is an isomorphism. Hence part (2) of the theorem is proven and part (1) follows from part (2) and Proposition~\ref{surj}.
\endproof
\end{thm}
\begin{example}
For $n=1$ and $k=2$, the fusion product $V(m_1)*V(m_2)$ can be filtered by
$$V_{\mathfrak{sl}_{2}}(m_1)*V_{\mathfrak{sl}_{2}}(m_2)\oplus V_{\mathfrak{sl}_{2}}(m_1-1)*V_{\mathfrak{sl}_{2}}(m_2-1)\oplus V_{\mathfrak{sl}_{2}}(m_1)*V_{\mathfrak{sl}_{2}}(m_2-1)\oplus V_{\mathfrak{sl}_{2}}(m_1-1)*V_{\mathfrak{sl}_{2}}(m_2).$$
\end{example}
\subsection{}
The following corollary is immediate.
\begin{cor}\mbox{}
\begin{enumerate}
\item The fusion product $V(m_1\delta_1)*\cdots*V(m_k\delta_1)$ is independent of the parameters.
\item Let $\blambda=(m_1\delta_1,\dots,m_k\delta_1), \bmu=(n_1\delta_1,\dots,n_k\delta_1)\in P^+(|\mathbf{m}|\delta_1,k)$ such that $\blambda\preceq \bmu$. Then we have a surjective map of $\mathbf{U}(\Lg[t])$--modules
$$V(n_1\delta_1)*\cdots*V(n_k\delta_1)\twoheadrightarrow V(m_1\delta_1)*\cdots*V(m_k\delta_1).$$
\item The $q$--character of $V(m_1\delta_1)*\cdots*V(m_k\delta_1)$ is given by 
$$\sum_{0\leq \ell \leq k}\ \sum_{0\leq i_1 <\cdots <i_k\leq k-1} q^{i_1+\cdots+i_{\ell}}\ \ch_q \Big(V_{\Lg_0}\big(\varphi_{i_1}\circ \cdots \circ \varphi_{i_\ell}(m_1)\delta_1\big)*\cdots *V_{\Lg_0}\big(\varphi_{i_1}\circ \cdots \circ \varphi_{i_\ell}(m_k)\delta_1\big)\Big).$$
\end{enumerate}
\end{cor}

For the Lie superalgebra $\mathfrak{osp}(1,2)$ we can give a PBW type basis. From Theorem~\ref{mainthm} and \cite[Theorem 5]{CV13} we get that the following set forms a basis of the fusion product

$$\bigcup_{0\leq \ell\leq k}\ \bigcup_{0\leq i_1 <\cdots <i_\ell\leq k-1}\Big\{(x_{2\delta_1}^{-}\otimes 1)^{j_1}\cdots (x_{2\delta_1}^{-}\otimes t^{k-1})^{j_k}(x_{\delta_1}^{-}\otimes t^{i_1})\cdots (x_{\delta_1}^{-}\otimes t^{i_{\ell}})\mid (j_1,\dots,j_k)\in S(\mathbf m, i_1,\dots,i_\ell)\Big\},$$
where $(j_1,\dots,j_k)\in S(\mathbf m, i_1,\dots,i_\ell)$ if and only if for all $2\leq r \leq k+1$ and $1\leq s\leq r-1$

$$(sj_{r-1}+(s+1)j_r)+2\sum_{p=r+1}^sj_p\leq \sum_{p=r-s}^{k}\varphi_{i_1}\circ\cdots\circ \varphi_{i_\ell}(m_p).$$
%
\section{Truncated Weyl modules and Demazure type modules}\label{section5}
In this section we discuss the connection of fusion products with truncated Weyl modules for $\Lg=\mathfrak{osp}(1,2)$. For the connection of theses modules for finite--dimensional simple Lie algebras we refer to \cite{KL14}. For the rest of this section we set $\alpha=2\delta_1$.
\subsection{}
\begin{defn}
For $n,N\in \mathbb{Z}_+$, the truncated Weyl module $\W(n,N)$ is the $\Lg[t]$--module generated by a non--zero element $w_{n,N}$ with defining relations:
\begin{align*}&\Ln^+[t]w_{n,N}=0,\ (h_{\alpha}\otimes t^r)w_n=n\delta_{r,0}w_{n,N},\ \forall r\in \mathbb{Z}_+,&\\& (x_{\alpha}^{-}\otimes 1)^{n+1}w_{n,N}=0,\ (\Lg\otimes t^N\mathbb{C}[t])w_{n,N}=0.\end{align*}
\end{defn}

The truncated Weyl module can be defined for any Lie superalgebra $\Lg$. For a finite--dimensional simple Lie algebra $\mathfrak{t}$ it has been proven in certain cases (for details see \cite{KL14}) that the truncated Weyl module can be realized as a fusion product of finite--dimensional irreducible $\mathfrak{t}$--modules. Similar results hold for the Lie superalgebra $\mathfrak{osp}(1,2)$.

\begin{prop}\label{fustr}
Let $n\in \mathbb Z_{+}$ and write $n=kN+j$ for $0\leq j <N$. 
We have an isomorphism of $\mathbf U(\Lg[t])$--modules
$$\W(n,N)\cong V(k)^{*(N-j)}*V(k+1)^{*j}.$$ 
\proof
By using the presentation of fusion products in Theorem~\ref{mainthm} the proof is exactly the same as the one in \cite[Section 4.3]{KL14}.
\endproof
\end{prop}
In particular, we obtain for all $N\geq n$ that the truncated Weyl module $\W(n,N)$ is isomorphic to the local Weyl module $\W(n)$. The connection of $\W(n)$ with fusion products was already studied in \cite{FM15}.
\subsection{}

For a simple finite--dimensional Lie algebra $\mathfrak{t}$ the Demazure module $\D_{\mathfrak{l}}(\ell,\lambda)$ of level $\ell$ and highest weight $\lambda$, where $\lambda$ is a dominant integral weight for $\mathfrak{t}$, is a finite--dimensional $\mathfrak{t}[t]$--submodule of an integrable level $\ell$ representation for the corresponding affine Kac--Moody algebra $\widehat{\mathfrak{t}}$. Demazure modules can be presented as cyclic modules that have an explicit description of the annihilator of the generating element. We will recall the presentation only for $\mathfrak{sl}_2$. The Demazure module $\D_{\mathfrak{sl}_2}(\ell, n)$ is the $\mathfrak{sl}_{2}[t]$--module generated by a non--zero vector $v_{\ell,n}$ subject to the defining relations
\begin{align*}&\Ln_0^+[t]v_{\ell,n}=0,\ (h_{\alpha}\otimes t^r)v_{\ell,n}=n\delta_{r,0}v_{\ell,n},\ (x_{\alpha}^{-}\otimes t^r)^{\max\{0,n-\ell r\}+1}v_{\ell,n}=0, \ \forall r\in \mathbb{Z}_+.\end{align*}

We can define the analogue modules for the Lie superalgebra $\Lg[t]$ and call them Demazure type modules, namely 
$\D(\ell, n)$ is the $\Lg[t]$--module generated by a non--zero vector $v_{\ell,n}$ subject to the defining relations
\begin{align*}&\Ln^+[t]v_{\ell,n}=0,\ (h_{\alpha}\otimes t^r)v_{\ell,n}=n\delta_{r,0}v_{\ell,n},\ (x_{\alpha}^{-}\otimes t^r)^{\max\{0,n-\ell r\}+1}v_{\ell,n}=0, \ \forall r\in \mathbb{Z}_+.\end{align*}
We can prove the following; similar results are true for finite--dimensional simple Lie algebras (see \cite{CP01,CV13}).
\begin{prop}
Write $n=(q-1)\ell+m$ for some $0<m \leq \ell$. We have an isomorphism of $\mathbf{U}(\Lg[t])$--modules 

$$\D(\ell,n)\cong V(\ell)^{*(q-1)}*V(m)\ \mbox{ and }\ \W(n)\cong \D(1,n).$$
\proof
The second isomorphism follows from Proposition~\ref{fustr} and the first isomorphism. The first isomorphism can be easily deduced from \cite[Theorem 1]{CV13}, but we will present the proof for the convenience of the reader. The highest weight of the fusion product satisfies obviously the defining relations of $D(\ell,n)$ and hence we have a surjective map 
$D(\ell,n)\twoheadrightarrow V(\ell)^{*(q-1)}*V(m)$. In order to prove a surjective map into the other direction let $(b_k)_{k\geq 0}\in _p\mathbf{S}(r,s)$ where $r,s\in \mathbb{N}$, $p\in \mathbb{Z}_+$ such that 
$$r+s\geq \begin{cases}
1+rp+(q-p-1)\ell+m, & \text{ if $p\leq q-1$}\\
1+rp,& \text{ else.}\end{cases}
$$
Since $(x_{\alpha}^{-}\otimes t^q)v_{\ell,n}=0$ we can assume that $b_k=0$ for all $k\geq q$. It follows that $p\leq q-1$, because otherwise we get a contradiction
$$qr\geq r+s \geq 1+rp \geq 1+rq.$$  
In this case we get 
$$r+s\geq 1+rp+(q-p-1)\ell+m\geq \begin{cases} 1+(q-1)r+m,& \text{ if $\ell \geq r$}\\ 1+n,& \text{ else.} \end{cases}$$
If $r>\ell$, the statement is obvious; so let $\ell \geq r$. Now the proposition follows with
$$1+(q-1)r+m\leq r+s \leq r+(q-1)b_{q-1}+(r-b_{q-1})(q-2)=b_{q-1}+r(q-1) \Rightarrow  m+1\leq b_{q-1}.$$
\end{prop}

%
\section{Proof of \texorpdfstring{Theorem~\ref{appe}}{Theorem}}\label{section6}
This section is dedicated to the proof of Theorem~\ref{appe}. Recall that $\Lg=\Lg_0\oplus \Lg_1$ and $\Lg_0$ is isomorphic to the symplectic Lie algebra $\mathfrak{sp}_{2n}$.
\subsection{}
We shall prove that we have an isomorphism of $\mathbf{U}(\Lg_0[t])$--modules 
\begin{equation}V_{\Lg_0}(m_1\delta_1)*\cdots * V_{\Lg_0}(m_k\delta_1)\cong \W_{\Lg_0}(|\mathbf{m}|\delta_1)/\mathfrak{I}(\mathbf{m}).\end{equation}
Recall that $\mathfrak{I}(\mathbf{m})$ is the submodule generated by the elements
\begin{equation}\label{uuu}\mathbf{x}^-_{\alpha}(r,s),\ \mbox{ for all }\alpha\in R^+_0,\ r,s\in \mathbb{N} : \ s+r\ge 1+ rp+\sum_{j\ge p+1}m_j\delta_1(h_{\alpha}),\ \mbox{ for some $p\in \mathbb{Z}_+$}.\end{equation}
We denote by $R_0^+(\neq 0)$ the set of all positive even roots $\beta$ such that $\delta_1(h_{\beta})\neq 0$ and by $R_0^+(0)$ the complement of $R_0^+(\neq 0)$ in $R^+_0$.
Further, for $\boldsymbol \nu=(\nu_{\alpha})_{\alpha\in R_0^+(\neq 0)}$ and $\ell\in \mathbb Z_+$ define 
$$\mu(\boldsymbol \nu,\mathbf m)=\sum_{j=1}^k\max\{0,|\boldsymbol\nu|-m_j\},\quad \mathbf S_{\boldsymbol \nu}(\ell)=\big\{\mathbf j= (j_{\alpha})_{\alpha\in R_0^+(\neq 0)}\mid j_{\alpha}\leq (k-1)\nu_{\alpha}, |\mathbf{j}|=(k-1)|\boldsymbol\nu|-\ell\big\}.$$
\subsection{}
\begin{proof}
It is easy to see that we have a surjective homomorphism 
$$\W_{\Lg_0}(|\mathbf{m}|\delta_1)\twoheadrightarrow V_{\Lg_0}(m_1\delta_1)*\cdots * V_{\Lg_0}(m_k\delta_1).$$
Let $\alpha\in R_0^+$ and consider the $\mathfrak{sl}_2$--triple $\{x^{\pm}_{\alpha}, h_{\alpha}\}$ and the corresponding current algebra $\mathfrak{sl}_{2}[t]$. We have a surjective map of $\mathbf{U}(\mathfrak{sl}_{2}[t])$--modules 
\begin{equation}\label{sss}V_{\mathfrak{sl}_{2}}(m_1\delta_1(h_{\alpha}))*\cdots * V_{\mathfrak{sl}_{2}}(m_k\delta_1(h_{\alpha}))\twoheadrightarrow \mathbf{U}(\mathfrak{sl}_{2}[t])(v_{m_1\delta_1}*\cdots*v_{m_k\delta_1}),\end{equation}
where the right hand side of \eqref{sss} is considered as a subspace of $V_{\Lg_0}(m_1\delta_1)*\cdots * V_{\Lg_0}(m_k\delta_1)$. It has been proved in \cite[Section 6]{CV13} that all elements in \eqref{uuu} act by zero on the highest weight vector of the fusion product $V_{\mathfrak{sl}_{2}}(m_1\delta_1(h_{\alpha}))*\cdots * V_{\mathfrak{sl}_{2}}(m_k\delta_1(h_{\alpha}))$ and hence by \eqref{sss} also on $v_{m_1\delta_1}*\cdots*v_{m_k\delta_1}$. It follows that we have a surjective map 
$$\W_{\Lg_0}(|\mathbf{m}|\delta_1)/\mathfrak{I}(\mathbf{m})\twoheadrightarrow V_{\Lg_0}(m_1\delta_1)*\cdots * V_{\Lg_0}(m_k\delta_1).$$
In order to prove the theorem, it remains to verify that the dimensions coincide. We consider a filtration on $V_{\Lg_0}(m\delta_1)=\mathbf{U}(\Ln_0^-)v_{m\delta_1}$ defined by
  \begin{equation*}\label{PBWfiltration}\mathbf{U}(\Ln_0^-)_sv_{m\delta_1}, \mbox{ where }\mathbf{U}(\Ln_0^-)_s=\spa\{x_{1}\cdots x_{l}\mid x_j\in\Ln_0^-, l\leq s\},\end{equation*}
and build the associated graded space $\gr V_{\Lg_0}(m\delta_1)$ with respect to this filtration. From \cite[Corollary 3.8]{FFL2011} and a straightforward calculation it follows that we have an isomorphism of $S(\Ln_0^{-})$--modules 
$$\gr V_{\Lg_0}(m\delta_1)\cong S(\Ln^{-})/\mathcal{I}(m\delta_1),$$
where $S(\Ln_0^{-})$ denotes the symmetric algebra of $\Ln_0^{-}$ and $\mathcal{I}(m\delta_1)$ the ideal generated by the elements 
$$\Big\{x_{\beta}^-,\ \prod_{\gamma\in R_0^+(\neq 0)}(x_{\gamma}^-)^{s_{\gamma}}\mid \mbox{ for all } \beta\in R_0^+(0) \mbox{ and } \mathbf{s}=(s_{\gamma})_{\gamma\in R_0^+(\neq 0)}\in \mathbb{Z}_+^{2n-1} \mbox{ such that $\sum_{\gamma}s_{\gamma}=m+1$}\Big\}.$$
In the language of \cite[Section 2.2]{FJKLM04} we have shown that $\gr V_{\Lg_0}(m\delta_1)\cong V_{m}^{(2n)}$ and hence we obtain with \cite[Theorem 2.4]{FJKLM04} that
$$\gr V_{\Lg_0}(m_1\delta_1)*\cdots *\gr V_{\Lg_0}(m_k\delta_1)\cong S(\Ln_0^-[t])/\mathcal{J}(|\mathbf{m}|),$$
where $\mathcal{J}(|\mathbf{m}|)$ is the ideal generated by the elements 
\begin{equation}\label{glach55}\big\{x^{-}_{\alpha}\otimes t^r\mid r\geq k,\ \alpha\in R_0^+\big\},\ \big\{x^{-}_{\alpha}\otimes t^r\mid r\geq 0,\ \alpha\in R_0^+(0)\big\}\end{equation} and
\begin{equation}\label{gleichung00}
\sum_{\mathbf{j}\in\mathbf S_{\boldsymbol \nu}(\ell)}\sum_{(b^{\alpha}_{k})_{k\geq 0}\in \mathbf{S}(\nu_\alpha,j_\alpha)}\prod_{\alpha\in R_0^+(\neq 0)}(x^{-}_{\alpha}\otimes 1)^{(b^{\alpha}_0)}\cdots (x^{-}_{\alpha}\otimes t^{|\mathbf j|})^{(b^{\alpha}_{|\mathbf j|})},\end{equation}
for all $\boldsymbol\nu=(\nu_{\alpha})_{\alpha\in R_0^+(\neq 0)}$ and $\ell <\mu(\boldsymbol \nu,\mathbf m)$. In the remaining part of the proof we will show that the generators of $\mathcal{J}(|\mathbf{m}|)$ act by zero on the highest weight vector of some PBW graded version of $\W_{\Lg_0}(|\mathbf{m}|\delta_1)/\mathfrak{I}(\mathbf{m})$, which would immediately imply $\dim \W_{\Lg_0}(|\mathbf{m}|\delta_1)/\mathfrak{I}(\mathbf{m})\leq \prod_{j=1}^k\dim V_{\Lg_0}(m_j\delta_1)$ and hence the theorem. 
Let $\mathbf{U}(\Ln_0^-[t])_{i}$ be the span of all monomials of the form
$$
\big(x_1\otimes p_1(t)\big)\cdots \big(x_r\otimes p_r(t)\big)\ 
\text{such that\ }r\le i.
$$
and consider the induced increasing filtration on $\W_{\Lg_0}(|\mathbf{m}|\delta_1)/\mathfrak{I}(\mathbf{m})$:
$$\cdots \subset \mathbf U(\Ln_0^-[t])_{s}w_{|\mathbf{m}|\delta_1}\subset \mathbf U(\Ln_0^-[t])_{s+1}w_{|\mathbf{m}|\delta_1} \subset $$

We build the associated graded space with respect to the above filtration and note that it is a cyclic $S(\mathfrak{n}_0^{-}[t])$--module whose cyclic vector $\gr(w_{|\mathbf{m}|\delta_1})$ is annihilated by the elements \eqref{glach55}. Further we have an $\mathbf U(\mathfrak{n}_0^+[t])$ action on the associated graded space which is induced by the $\mathbf U(\mathfrak{n}_0^+[t])$ action on $\W_{\Lg_0}(|\mathbf{m}|\delta_1)/\mathfrak{I}(\mathbf{m})$. Let $\boldsymbol \nu =(\nu_{\alpha})_{\alpha\in R^+_0(\neq 0)}$ and $\ell\in \mathbb Z_+$ and $\mathbf j\in \mathbf S_{\boldsymbol \nu}(\ell)$. Our aim is to prove that the elements \eqref{gleichung00} act by zero on $\gr(w_{|\mathbf{m}|\delta_1})$. Let $p\in \mathbb{Z}_+$ minimal such that $|\boldsymbol \nu|\geq m_{p+1}$. Since
$$|\boldsymbol \nu|+|\mathbf j|\geq k|\boldsymbol \nu|-\ell\geq 1+k|\boldsymbol \nu|-\mu(\boldsymbol \nu,\mathbf m)\geq 1+p|\boldsymbol \nu|+\sum_{j\geq p+1} m_{j}$$
we obtain for all $\alpha\in R^+_0(\neq 0)$
$$\mathbf U(\Ln^+[t])\circ \mathbf x_{\alpha}^-\big(|\boldsymbol \nu|,|\mathbf j|\big)\gr(w_{|\mathbf{m}|\delta_1})=0.$$
But an easy calculation shows that
$$\prod_{\substack{\beta\in R_0^+(\neq 0)\\ \beta\neq 2\delta_1}} (x_{\beta}^{+}\otimes 1)^{(\nu_{2\delta_1-\beta})}\circ \mathbf{x}_{2\delta_1}^-\big(|\boldsymbol \nu|,|\mathbf j|\big)\gr(w_{|\mathbf{m}|\delta_1})$$
gives the element \eqref{gleichung00}.
\end{proof}

%

\bibliographystyle{plain}
\bibliography{liesuperfus}

\begin{thebibliography}{10}

\bibitem{BCSV15}
Rekha Biswal, Vyjayanthi Chari, Lisa Schneider, and Sankaran Viswanath.
\newblock Demazure {F}lags, {C}hebyshev polynomials, {P}artial and {M}ock theta
  functions.
\newblock arXiv:1502.05322.

\bibitem{CLS15}
Lucas Calixto, Joel Lemay, and Savage.
\newblock Weyl modules for {L}ie superalgebras.
\newblock arXiv:1505.06949.

\bibitem{CFS12}
Vyjayanthi Chari, Ghislain Fourier, and Daisuke Sagaki.
\newblock Posets, tensor products and {S}chur positivity.
\newblock {\em Algebra Number Theory}, 8(4):933--961, 2014.

\bibitem{CIK14}
Vyjayanthi Chari, Bogdan Ion, and Deniz Kus.
\newblock Weyl modules for the hyperspecial current algebra.
\newblock {\em Int. Math. Res. Not. IMRN}, (15):6470--6515, 2015.

\bibitem{CL06}
Vyjayanthi Chari and Sergei Loktev.
\newblock Weyl, {D}emazure and fusion modules for the current algebra of
  {$\mathfrak{sl}_{r+1}$}.
\newblock {\em Adv. Math.}, 207(2):928--960, 2006.

\bibitem{CP01}
Vyjayanthi Chari and Andrew Pressley.
\newblock Weyl modules for classical and quantum affine algebras.
\newblock {\em Represent. Theory}, 5:191--223 (electronic), 2001.

\bibitem{CSSW14}
Vyjayanthi Chari, Lisa Schneider, Perri Shereen, and Jeffrey Wand.
\newblock Modules with {D}emazure {F}lags and {C}haracter {F}ormulae.
\newblock {\em SIGMA Symmetry Integrability Geom. Methods Appl.}, 10, 2014.

\bibitem{CSVW14}
Vyjayanthi Chari, Perri Shereen, R.~Venkatesh, and Jeffrey Wand.
\newblock A {S}teinberg type decomposition theorem for higher level {D}emazure
  modules.
\newblock arXiv:1408.4090.

\bibitem{CV13}
Vyjayanthi Chari and R.~Venkatesh.
\newblock Demazure modules, {F}usion products and {Q}--systems.
\newblock arXiv:1305.2523 , Communications in {M}athematical {P}hysics, to
  appear.

\bibitem{FF02}
B.~Feigin and E.~Feigin.
\newblock {$Q$}-characters of the tensor products in {$\mathfrak{sl}_2$}-case.
\newblock {\em Mosc. Math. J.}, 2(3):567--588, 2002.
\newblock Dedicated to Yuri I. Manin on the occasion of his 65th birthday.

\bibitem{FJKLM04}
B.~Feigin, M.~Jimbo, R.~Kedem, S.~Loktev, and T.~Miwa.
\newblock Spaces of coinvariants and fusion product. {II}.
  {$\widehat{\mathfrak{sl}}\sb 2$} character formulas in terms of {K}ostka
  polynomials.
\newblock {\em J. Algebra}, 279(1):147--179, 2004.

\bibitem{FL99}
Boris Feigin and Sergei Loktev.
\newblock On generalized {K}ostka polynomials and the quantum {V}erlinde rule.
\newblock In {\em Differential topology, infinite-dimensional {L}ie algebras,
  and applications}, volume 194 of {\em Amer. Math. Soc. Transl. Ser. 2}, pages
  61--79. Amer. Math. Soc., Providence, RI, 1999.

\bibitem{FFL2011}
Evgeny Feigin, Ghislain Fourier, and Peter Littelmann.
\newblock P{BW} filtration and bases for symplectic {L}ie algebras.
\newblock {\em Int. Math. Res. Not. IMRN}, (24):5760--5784, 2011.

\bibitem{FM15}
Evgeny Feigin and Ievgen Makedonskyi.
\newblock {W}eyl modules for $\mathfrak{osp}(1,2)$ and nonsymmetric {M}acdonald
  polynomials.
\newblock arXiv:1507.01362.

\bibitem{FG11}
R.~Fioresi and F.~Gavarini.
\newblock On the construction of {C}hevalley supergroups.
\newblock In {\em Supersymmetry in mathematics and physics}, volume 2027 of
  {\em Lecture Notes in Math.}, pages 101--123. Springer, Heidelberg, 2011.

\bibitem{FK11}
Ghislain Fourier and Deniz Kus.
\newblock Demazure modules and {W}eyl modules: {T}he twisted current case.
\newblock {\em Trans. Amer. Math. Soc.}, 365(11):6037--6064, 2013.

\bibitem{FoL07}
Ghislain Fourier and Peter Littelmann.
\newblock Weyl modules, {D}emazure modules, {KR}-modules, crystals, fusion
  products and limit constructions.
\newblock {\em Adv. Math.}, 211(2):566--593, 2007.

\bibitem{IK01}
Kenji Iohara and Yoshiyuki Koga.
\newblock Central extensions of {L}ie superalgebras.
\newblock {\em Comment. Math. Helv.}, 76(1):110--154, 2001.

\bibitem{K78}
V.~Kac.
\newblock Representations of classical {L}ie superalgebras.
\newblock In {\em Differential geometrical methods in mathematical physics,
  {II} ({P}roc. {C}onf., {U}niv. {B}onn, {B}onn, 1977)}, volume 676 of {\em
  Lecture Notes in Math.}, pages 597--626. Springer, Berlin, 1978.

\bibitem{KL14}
Deniz Kus and Peter Littelmann.
\newblock Fusion products and toroidal algebras.
\newblock arXiv:1411.5272 (Pacific Journal of Mathematics, to appear).

\bibitem{KV14}
Deniz Kus and R.~Venkatesh.
\newblock Twisted {D}emazure modules, {F}usion product decomposition and
  twisted {Q}--systems.
\newblock arXiv:1409.1201.

\bibitem{Mu12}
Ian~M. Musson.
\newblock {\em Lie superalgebras and enveloping algebras}, volume 131 of {\em
  Graduate Studies in Mathematics}.
\newblock American Mathematical Society, Providence, RI, 2012.

\bibitem{Na11}
Katsuyuki Naoi.
\newblock Weyl modules, {D}emazure modules and finite crystals for non-simply
  laced type.
\newblock {\em Adv. Math.}, 229(2):875--934, 2012.

\end{thebibliography}
\end{document}